\documentclass[12pt]{smfart}
\usepackage[frenchb]{babel}
\usepackage{mathptmx}
\usepackage{amsmath}
\usepackage{amssymb}
\usepackage{smfthm}
\usepackage{amscd}
\usepackage{smfenum}

\usepackage{amssymb,amsfonts}

\usepackage[T1]{fontenc}     
\usepackage[latin1]{inputenc}

\newcounter{ctrmasection}
\newcounter{ctrmasoussection}

\newenvironment{pf*}{  \bf{D{\'e}monstration : }\rm}{$\Box$}

{\begin{displaymath} \begin{array}{c}}%
{\end{array} \end{displaymath}}
\newcommand{\Aut}{\operatorname{Aut}}

\newcommand{\Spec}{\operatorname{Spec}}
\newcommand{\GTHS}{\operatorname{GTHS}}
\newcommand{\GTRS}{\operatorname{GTRS}}

\begin{document}
\title[Dynamique $p$-adique arithm\'etique des automorphismes affines]{Sur la dynamique  $p$-adique arithm\'etique des automorphismes de l'espace affine} 
\alttitle{On the  arithmetic $p$-adic dynamics of automorphisms of the affine space}
\author{Sandra Marcello}

\address{Max-Planck-Institut f\"ur Mathematik,
Vivatsgasse 7, 53111, Bonn, Deutschland \\}
\email{marcello@mpim-bonn.mpg.de}
\date{\today}

\begin{abstract}
Nous montrons que l'ensemble des p\'eriodes d'un automorphisme d'un sous-groupe du groupe des automorphismes du plan affine d\'efini sur un corps $p$-adique est major\'e par une constante ind\'ependante de l'automorphisme. Nous en d\'eduisons  que ce sous groupe v\'erifie la conjecture de dynamique arithm\'etique de Silverman, ce r\'esultat est inclus dans un r\'esultat de L. Denis \cite{denis}. 
\end{abstract}

\begin{altabstract}
We show that the set of periods for a automorphism of a sub-group of the group of automorphisms of the affine plane defined over a $p$-adic field is bounded above by a constant independent from the automorphism. We deduce from this result, that the Silverman's conjecture on arithmetic dynamics is true for this subgroup, this result is due to L. Denis \cite{denis}. 
\end{altabstract}


\altkeywords{Iterates, automorphisms of affine space, $p$-adic fields}

\maketitle


\mainmatter

\section{Abridged English version}

We study, over a $p$-adic field (a finite extension of $\mathbb{Q}_p$), the periodic points of automorphisms of affine spaces.

\paragraph*{Notations and definitions}
\begin{itemize}
\item $p$ a prime  number, $K$ a $p$-adic field, $\mathcal{O}$ the ring of integers of $K$,
\item $\Aut(\mathbb{A}^r)(K)$ the group of  automorphisms of affine space of dimension $r\geq 2$ defined over $K$. 
\item Let $\phi\in \Aut(\mathbb{A}^r)(K)$. We also write $\phi$  for tha rational map of the projective scheme associated
to  ${\mathbb{A}}^r$. We write $Z(\phi)$ the locus of non-definition of the rational map associated to $\phi$.
Let $P$ be a $\phi$-periodic point. We write $n_{P,\phi}$ the {\it{$\phi$-period}} of $P$ {\it{i.e}} the smallest integer  $n\in \mathbb{N}\setminus\{0\}$ such that $\phi^n(P)=P$.
\item A {\it generalised  H{\'e}non map} is a map $g$ of the type $g(x,y)=(p(x)-ay,x)$ where $a$ is a non-zero constant and $p$ is a monic 
 polynomial of  degree at least 2.
\item An automorphism $\phi$ of $\mathbb{A}^r(K)$ is said {\it{triangular}}
  if it exists
$F_i\in K[X_{i+1},\dots,X_{r}]$ for $1\leq i\leq r-1$, $F_r\in K$  and  $a_i\in
K^*$ for $1\leq i\leq r$ such that\,:
\begin{eqnarray*}
\phi(X_{1},\dots,X_{r})&=&(a_1X_{1}+F_1(X_{2},\dots,X_{r}),a_2X_2+F_2(X_{3},\dots,X_{r}),\\
&               &\quad \dots,a_{r-1}X_{r-1}+F_{r-1}(X_{r}),a_rX_r+F_r).
\end{eqnarray*}
\item  Let  $\mathcal{X}$ be an integer proper scheme of finite type fini over $\Spec( \mathcal{O})$, we write $\mathcal{X}_s$ for the special fiber of  $\mathcal{X}$.
\end{itemize}

\paragraph*{Definition}
Let $\psi\in \Aut(\mathbb{A}^2)(K)$. We say that $\psi$ is a {\it{ special H\'enon map}} if : 
\begin{itemize}
\item $\psi$ is a finite product of generalised  H\'enon maps,
\item $$Z(\psi|_{\mathcal{X}_s})\cap Z(\psi^{-1}|_{\mathcal{X}_s})=\emptyset.$$
\end{itemize}
We define the group $GTHS(K)\subset \Aut(\mathbb{A}^2)(K)$ :
\begin{eqnarray*}
\GTHS(K)=&\{& f^{-1}\phi f \quad \mbox{avec}\quad \phi,f\in \Aut(\mathbb{A}^2)(K)\quad 
 \mbox{and $\phi$  triangular or special H\'enon map}
\}.
\end{eqnarray*}

Our main result is:

\paragraph*{{\bf{Theorem  A}}}
{\it{Let $K$ be a $p$-adic field. It exists $M(K)\in \mathbb{N}$ such that  $\forall \phi\in \GTHS(K)$ and for all   $\phi$-p\'eriodic point  $P$ we have :}}
$$0\leq n_{P,\phi} \leq M.$$

\section{Introduction}
 Apr\`es J. Silverman \cite{Silv1} et L. Denis \cite{denis}, dans \cite{1art} et \cite{2art} nous nous sommes int\'eress\'ee \`a l'\'etude des propri\'et\'es des it\'er\'es des automorphismes de l'espace affine  d\'efinis sur un corps de nombres.  
 En 1995, P. Morton et J. Silverman se sont notamment int\'eress\'es au comportement apr\`es r\'eduction modulo $p$ des points p\'eriodiques d'applications rationnelles de la droite projective \cite{ms1}, de plus, 
ces derni\`eres ann\'ees l'\'etude de la dynamique ultram\'etrique ou $p$-adique (sur $\mathbb{C}_{p}$) a pris de l'essor (voir par exemple \cite{bened1},\cite{bened2} et \cite{rl1}).
Ainsi l'\'etude du comportement des it\'er\'es des automorphismes de l'espaces affines sur un corps $p$-adique appara\^\i t comme naturelle.
La notion de point p\'eriodique peut \^etre vue comme l'{\it{analogue}} de la notion de point de torsion des vari\'et\'es ab\'eliennes. Les r\'esultats de finitude de points de torsion de vari\'et\'es ab\'eliennes que ce soit sur un corps nombres (voir par exemple \cite{Hindry} partie C) ou sur un corps p-adique (\cite{Milne} remarque 8.4) sont classiques.
Nous nous int\'eressons ici \`a ce m\^eme type de questions dans le cas d'un corps $p$-adique, nous obtenons m\^eme mieux \`a savoir une {\it{majoration uniforme}}.

\paragraph*{Notations et d\'efinitions}
Nous noterons:
\begin{itemize}
\item $p$ un nombre premier, $K$ un corps $p$-adique (\`a savoir une extension finie de $\mathbb{Q}_{p}$), $\mathcal{O}$ l'anneau des entiers de $K$,
\item $\Aut(\mathbb{A}^r)(K)$ le groupe des automorphismes de l'espace affine de dimension $r\geq 2$ d\'efini sur $K$. 
\item Soit $\phi\in \Aut(\mathbb{A}^r)(K)$. Par abus de langage, nous noterons \'egalement $\phi$ l'application birationnelle du sch\'ema projectif associ\'e \`a  ${\mathbb{A}}^r$. Nous d{\'e}signons par $Z(\phi)$ le lieu de non-d{\'e}finition de l'application rationnelle associ{\'e}e {\`a} $\phi$.
Soit $P$ un point $\phi$-p\'eriodique. Nous noterons $n_{P,\phi}$ la {\it{$\phi$-p\'eriode}} de $P$ {\it{i.e}} le plus petit entier $n\in \mathbb{N}\setminus\{0\}$ tel que $\phi^n(P)=P$.
\item Une {\it application de H{\'e}non g{\'e}n{\'e}ralis{\'e}e} est une application $g$ de la forme $g(x,y)=(p(x)-ay,x)$ o{\`u} $a$ est une constante non nulle et $p$ est un 
 polyn\^ome unitaire de degr{\'e} au moins 2.
\item Application  {\it{triangulaire}} voir la d\'efinition \ref{deftri}.
\item  Soit $\mathcal{X}$ un sch\'ema int\`egre propre de type fini sur $\Spec( \mathcal{O})$, nous noterons $\mathcal{X}_s$ la fibre sp\'eciale de  $\mathcal{X}$.
\end{itemize}

\begin{defi}\label{defgths}
Soit $\psi\in \Aut(\mathbb{A}^2)(K)$. Nous dirons que $\psi$ est une {\it{ application de H\'enon sp\'eciale}} si elle v\'erifie les conditions suivantes:
\begin{itemize}
\item $\psi$ est un produit fini d'applications de H\'enon g\'en\'eralis\'ees
\item $$Z(\psi|_{\mathcal{X}_s})\cap Z(\psi^{-1}|_{\mathcal{X}_s})=\emptyset.$$
\end{itemize}
Nous d\'efinissons le groupe $G(K)\subset \Aut(\mathbb{A}^2)(K)$, par:

$\GTHS(K)=\{ f^{-1}\phi f \quad \mbox{avec}\quad  
\phi,f\in \Aut(\mathbb{A}^2)(K)\quad $
\vspace{-0.3cm}
$$\mbox{et $\phi$ application triangulaire ou application de H\'enon sp\'eciale}
\}.$$ 
\end{defi}

Si on suppose que $\phi$ est une application triangulaire ou un produit fini d'application de H\'enon g\'en\'eralis\'ees, on obtient $\Aut(\mathbb{A}^2)(K)$. Ce r\'esultat a \'et\'e montr\'e par S. Friedland et J. Milnor \cite{FM} sur $\mathbb{C}$ et $\mathbb{R}$ \`a l'aide du th\'eor\`eme de Jung \cite{Jung}. Nous avons v\'erifi\'ee que la d\'emonstration de  S. Friedland et J. Milnor \cite{FM} est encore valable sur un corps $p$-adique, \`a l'aide du th\'eor\`eme de Jung qui a \'et\'e d\'emontr\'e  sur un corps quelconque par Makar-Limanov\cite{Makar}, une autre preuve (valable pour tout corps) de ce th\'eor\`eme se trouve dans\cite{Cohn}.

\paragraph*{{\bf{Th\'eor\`eme  A}}}
{\it{Soit $K$ un corps $p$-adique. Il existe $M(K)\in \mathbb{N}$ tel que  $\forall \phi\in \GTHS(K)$ et pour tout point $\phi$-p\'eriodique $P$ nous avons :}}
$$0\leq n_{P,\phi} \leq M.$$ 

 Nous obtenons gr\^ace \`a ce r\'esultat une forme faible d'un r\'esultat de dynamique arithm\'etique d\^u \`a L. Denis \cite{denis}.

\paragraph*{\textbf{D\'efinition}} Soit $F$ un corps de nombres. Soit $\phi\in\Aut(\mathbb{A}^r(F))$.

\noindent  Soit $P\in \mathbb{A}^r(F)$. Ce point $P$  
est un {\it{point $\phi$-p{\'e}riodique isol{\'e}}}
si et seulement s'il existe $n \in \mathbb{N}\setminus \{0\}$ tel que
$P$ est isol{\'e} (au sens de la topologie de Zariski) dans
 $$\{Q\in \mathbb{A}^r(F)\quad \mbox{tel que}\quad \phi ^{n}(Q)=Q\}.$$

\noindent{\textbf{Conjecture de Silverman:}} {\it{Soit $F$ un corps de nombres. Soit $\phi\in Aut(\mathbb{A}^r(F))$.
L'ensemble des points p{\'e}riodiques isol\'es de $\phi$ est un ensemble fini.}}

\paragraph*{{\bf{Corollaire B}}}
{\it{ La conjecture de Silverman est vraie pour $\GTHS(F)$.}}

L'approche utilis\'ee  pour d\'emontrer le th\'eor\`eme A est analogue \`a celle utilis\'ee par l'auteure  dans  \cite{1art}. Nous obtenons une majoration uniforme pour certains automorphismes r\'eguliers et les automorphismes triangulaires.
Nous obtenons \'egalement un analogue du th\'eor\`eme A en dimension sup\'erieure.

\noindent Dans le dernier paragraphe, nous faisons le lien avec la dynamique arithm\'etique et montrons le corollaire B.

\paragraph*{Remerciements}L'auteure remercie le r\'eseau europ\'een Arithmetic Algebraic Geometry.
Ce travail a \'et\'e effectu\'e au Max-Planck Institut f\"ur Mathematik que l'auteure remercie pour les excellentes conditions de travail. L'auteure remercie  tr\`es vivement Jaya Iyer de lui avoir parl\'e des travaux de N. Fakhruddin.
 L'auteure remercie N. Fakhruddin de lui avoir fait part  d'un probl\`eme  avec un lemme dans une version pr\'ec\'edente de ce texte.

\section{Automorphismes r\'eguliers}

\subsection{d\'efinition et propri\'et\'es}

La notion d'automorphisme r\'egulier a \'et\'e introduite par N. Sibony dans \cite{Sibony}, leur dynamique holomorphe \cite{Sibony} et arithm\'etique  \cite{0art} et \cite{1art} sont int\'eressantes.

\begin{defi}
Soit $\phi$ un automorphisme de $\mathbb{A}^r$, $Z(\phi)$ d{\'e}signe le lieu de non-d{\'e}finition de l'application
rationnelle (d{\'e}finie sur $\mathbb{P}^r$) associ{\'e}e {\`a} $\phi$. L'automorphisme $\phi$
est dit {\it{r{\'e}gulier}} 
si $\phi$ n'est pas une application lin\'eaire et \,:
$$Z(\phi)\cap Z(\phi^{-1})=\emptyset.$$
\end{defi}
Soit $H$ l'hyperplan {\`a} l'infini, nous avons $Z(\phi)\subset H$.

\begin{lemm}\label{lemmedenis}
Un produit fini d'applications de H{\'e}non g{\'e}n{\'e}ralis{\'e}es est un
automorphisme r{\'e}gulier.
\end{lemm}

\noindent Ce r{\'e}sultat s'obtient directement en composant les applications. 

 Nous utiliserons la proposition suivante que nous avons d\'emontr\'ee, dans le cas des corps de nombres, dans \cite{1art} (proposition 2.23), la d\'emonstration est \'egalement valable dans le cas des corps $p$-adiques ou des corps finis.

\begin{prop}\label{regn}
Soit $\phi$ un automorphisme r{\'e}gulier de ${\mathbb{A}}^r$.
 Alors, pour tout entier
naturel $n$ non nul $\phi^n$ est r{\'e}gulier et $Z(\phi^n)=Z(\phi)$.
\end{prop}

\begin{defi}
Soit $\phi$ un automorphisme r{\'e}gulier de ${\mathbb{A}}^r$. Nous dirons que $\phi$ est un {\it{ automorphisme r\'egulier sp\'ecial}} si :
$$Z(\psi|_{\mathcal{X}_s})\cap Z(\psi^{-1}|_{\mathcal{X}_s})=\emptyset.$$
\end{defi}

\subsection{Un r\'esultat pr\'eliminaire}
\begin{defi}
Soit $X$ une vari\'et\'e  sur un corps $K$. Soit $f:X\cdots \rightarrow X$ une application rationnelle. Nous dirons que $P$ est un point {\it{$f$-p\'eriodique}} si $f$ est d\'efinie en $P,f(P),f^2(P),\cdots$ et il existe $n>0$ tel que $f^n(P)=P$.
\end{defi}

Le th\'eor\`eme suivant est d\^u \`a N. Fakhruddin, nous donnons ici des \'etapes de la preuve et renvoyons le lecteur \`a \cite{fakh} pour plus de d\'etails.  

\begin{theo}\label{fakh}
Soit $\mathcal{X}$ un sch\'ema int\`egre propre de type fini sur $\Spec( \mathcal{O})$.
Alors, il existe une constante $M>0$ telle que 

\centerline{$\forall f:\mathcal{X} \cdots\rightarrow \mathcal{X}$ application birationnelle telle que $\forall n>0$, $Z(f^n)\cap Z(f^{-n})=\emptyset$}
\noindent $\forall P\in \mathcal{X}(K)$ $f$-p\'eriodique nous avons :
$$1\leq n_{f,P}\leq M$$
\end{theo}

Il suffit de montrer le r\'esultat pour une composante connexe de l'adh\'erence de Zariski de l'orbite d'un point p\'eriodique, en effet le corps r\'esiduel \'etant fini le nombre de composantes connexes de l'adh\'erence de l'orbite peut \^etre major\'e ind\'ependamment  du point et de l'application. \`A chaque composante connexe on associe de mani\`ere unique un point ferm\'e $p_i$ dont on majore la p\'eriode par une constante $n-i$ ne d\'ependant que du sch\'ema. On consid\`ere un point de l'orbite de d\'epart dont la sp\'ecialisation est ce point ferm\'e, on consid\`ere l'orbite de ce point sous l'action de $f^{n_i}$ et on majore la p\'eriode. 
La majoration s'obtient d'abord modulo une puissance de $p$, puis il s'agit de montrer que la puissance peut elle-m\^eme \^etre major\'ee.

\subsection{Un r\'esultat g\'en\'eral pour certains automorphismes r\'eguliers}

\begin{theo}\label{threg}Soit $K$ un corps $p$-adique.
\begin{itemize} 
\item Il existe $M\in \mathbb{N}$  tel que 
pour tout  une applicationde H\'enon sp\'eciale $\phi$, 
pour tout $P\in {\mathbb{A}}^r(K)$ $\phi$-p\'eriodique alors:
$$1\leq n_{\phi,P} \leq M.$$
\item Soit $r\geq 3$. Il existe $M\in \mathbb{N}$  tel que 
pour tout automorphisme r\'egulier sp\'ecial  $\phi\in {\mathbb{A}}^r(K) $,pour tout $P\in {\mathbb{A}}^r(K)$ $\phi$-p\'eriodique alors:
$$1\leq n_{\phi,P} \leq M.$$
\end{itemize}
\end{theo}

La preuve de ce th\'eor\`eme d\'ecoule du th\'eor\`eme \ref{fakh} et de la proposition \ref{regn}.

\section{Preuve du r\'esultat principal}

\subsection{Applications triangulaires}

\begin{defi}\label{deftri}
Un automorphisme $\phi$ de $\mathbb{A}^r(K)$ est dit {\it{triangulaire}}
  s'il existe
$F_i\in K[X_{i+1},\dots,X_{r}]$ pour $1\leq i\leq r-1$, $F_r\in K$  et  $a_i\in
K^*$ pour $1\leq i\leq r$ tel que\,:
\begin{eqnarray*}
\phi(X_{1},\dots,X_{r})&=&(a_1X_{1}+F_1(X_{2},\dots,X_{r}),a_2X_2+F_2(X_{3},\dots,X_{r}),\\
&               &\quad \dots,a_{r-1}X_{r-1}+F_{r-1}(X_{r}),a_rX_r+F_r).
\end{eqnarray*}
\end{defi}
\paragraph*{Notation} Soit $\mu_{K}$ le groupe (fini) des racines de l'unit\'e de $K$.

\begin{lemm}\label{tri} Il existe $M\in \mathbb{N}$ tel que pour  toute une application triangulaire  $\phi$ de $\mathbb{A}^r(K)$, pour tout point $\phi$-p\'eriodique $P$: $1\leq n_{\phi,P}\leq M$.
\end{lemm} 
Nous d\'eterminons les  it\'er\'es, et consid\'erons deux cas suivant que l'un au moins des coefficients $a_i$ est ou non une racine de l'unit\'e et concluons en utilisant  le fait que $\mu_K$ est un groupe fini.

Si dans la d\'efinition \ref{defgths}, nous nous pla{\c{c}}ons en dimension quelconque $r$ et consid\'erons des automorphismes r\'eguliers sp\'eciaux \`a la place des applications de H\'enon sp\'eciales, nous d\'efinissons alors le groupe: $GTRS(K,r)$.

\subsection{Conjugaison par un automorphisme affine}
Le lemme suivant  a {\'e}t{\'e} utilis{\'e} par L. Denis \cite{denis} dans le cas des corps de nombres, sa d\'emonstration est \'egalement imm\'ediate dans le cas des corps $p$-adiques.

\begin{lemm}\label{conj}
Soient $\psi,\phi$ des automorphismes de l'espace affine
${\mathbb{A}}^r(K)$. Supposons qu'il existe $f$ un automorphisme de
l'espace affine ${\mathbb{A}}^r(K)$ tel que $\psi=f^{-1}\phi f $ alors :

\centerline{$P$ est $\psi$-p{\'e}riodique si et seulement si $f(P)$ est $\phi$-p{\'e}riodique.}
\end{lemm}

\subsection{Preuve du th\'eor\`eme A}

La preuve du th\'eor\`eme A d\'ecoule imm\'ediatement  du th\'eor\`eme \ref{threg}, et des lemmes \ref{tri}, \ref{lemmedenis} et \ref{conj}.
 
De plus, nous obtenons un r\'esultat analogue en dimension quelconque.

\begin{theo}\label{theogtrs}
Soit $K$ un corps $p$-adique. Il existe $M(K)\in \mathbb{N}$ tel que  $\forall \phi\in \GTRS(K,r)$ et pour tout point $\phi$-p\'eriodique $P$ nous avons :
$$0\leq n_{P,\phi} \leq M.$$ 
\end{theo}

\section{Lien avec la dynamique arithm\'etique}

\begin{prop}\label{propdyn}Si  pour tout  corps $p$-adique $K$, il existe $M\in \mathbb{N}\setminus\{0\}$ tel que pour tout $\phi\in \Aut(\mathbb{A}^r(K))$, pour tout point $\phi$-p\'eriodique $P$, nous avons : $1\leq n_{\phi,P}\leq M$, alors,

si $F$ est un corps de nombres, pour tout $\phi\in Aut(\mathbb{A}^r(F))$:

\centerline{l'ensemble des points $\phi$-p\'eriodiques isol\'es est fini.} 
\end{prop}

Pour d\'emontrer ce r\'esultat on utilise la d\'ecomposition suivante de l 
'ensemble des points $\phi$-p\'eriodiques $B(\phi,K)=\cup_{N\geq 1}B_N(\phi,K)$, avec $B_N(\phi,K):=\{P\in\mathbb{A}^n(K); \phi^N(P)=P\}$, c'est un ferm\'e, il admet donc un nombre fini de points isol\'es. Le th\'eor\`eme A nous permet de dire que la r\'eunion est finie, de plus tout corps de nombres peut \^etre plong\'e dans un corps $p$-adique, d'o\`u le r\'esultat.

La proposition \label{propdyn} nous permet conjointement avec le th\'eor\`eme A d'obtenir le corollaire B. Nous obtenons \`a l'aide du th\'eor\`eme \ref{theogtrs}, l'analogue du corollaire B pour le groupe $\GTRS(K,r)$.

\end{document}